\documentclass[12pt]{article}

\usepackage{amsmath, epsfig, cite}
\usepackage{amssymb}
\usepackage{amsfonts}
\usepackage{latexsym}

\newtheorem{thm}{Theorem}[section]

\newtheorem{conj}[thm]{Conjecture}

\numberwithin{equation}{section}

\begin{document}


\begin{center}
{\Large\bf Proof of a congruence on sums of powers of\\[5pt] $q$-binomial coefficients}
\end{center}

\vskip 2mm \centerline{Victor J. W. Guo$^1$ and Ji-Cai Liu$^2$}
\begin{center}
{\footnotesize Department of Mathematics, Shanghai Key Laboratory of
PMMP, East China Normal University,\\ 500 Dongchuan Road, Shanghai
200241,
 People's Republic of China\\
$^1${\tt jwguo@math.ecnu.edu.cn,\quad
http://math.ecnu.edu.cn/\textasciitilde{jwguo}}  \quad {$^2$\tt jc2051@163.com} }
\end{center}


\vskip 0.7cm \noindent{\bf Abstract.} We prove that, if $m,n\geqslant 1$ and $a_1,\ldots,a_m$ are nonnegative integers, then
\begin{align*}
\frac{[a_1+\cdots+a_m+1]!}{[a_1]!\ldots[a_m]!}\sum^{n-1}_{h=0}q^h\prod_{i=1}^m{h\brack a_i}
\equiv 0\pmod{[n]},
\end{align*}
where $[n]=\frac{1-q^n}{1-q}$, $[n]!=[n][n-1]\cdots[1]$, and
${a\brack b}=\prod_{k=1}^b\frac{1-q^{a-k+1}}{1-q^k}$. The $a_1=\cdots=a_m$ case confirms a recent conjecture of Z.-W. Sun. We also show that,
if $p$ is a prime greater than $\max\{a,b\}$, then
\begin{align*}
\frac{[a+b+1]!}{[a]![b]!}\sum_{h=0}^{p-1}q^h{h\brack a}{h\brack b}
\equiv (-1)^{a-b} q^{ab-{a\choose 2}-{b\choose 2}}[p]\pmod{[p]^2}.
\end{align*}


\vskip 2mm
\noindent{\it MR Subject Classifications}: 11A07, 11B65, 05A10

\section{Introduction}
Recall that
the {\it $q$-binomial coefficients} ${n\brack k}$ (see \cite{GR}) are defined by
$$
{n\brack k}
=\begin{cases}
\displaystyle\frac{(1-q^n)(1-q^{n-1})\cdots (1-q^{n-k+1})}{(1-q)(1-q^2)\cdots(1-q^k)}, &\text{if $0\leqslant k\leqslant n$},\\[10pt]
0,&\text{otherwise.}
\end{cases}
$$
$q$-Binomial coefficients are closely related to binomial coefficients by the relation $\lim_{q\to 1}{n\brack k}={n\choose k}$.
Recently, Z.-W. Sun \cite{Sun}  proved many interesting congruences
on sums involving binomial coefficients or $q$-binomial coefficients.
For example, Sun \cite{Sun} proved that, for any nonnegative integers $n$ and $k$ with $n>k$, there holds
\begin{align}\label{eq:sun1}
[2k+1]{2k\brack k}\sum_{h=0}^{n-1}q^h {h\brack k}^2\equiv 0\pmod{[n]},
\end{align}
where $[n]:=1+q+\cdots+q^{n-1}$, and so
\begin{align}
(2k+1){2k\choose k}\sum_{h=0}^{n-1}{h\choose k}^2\equiv 0\pmod{n}. \label{eq:sun2}
\end{align}

He also made the following conjecture, which is a generalization of \eqref{eq:sun1} and \eqref{eq:sun2}.
\begin{conj}\label{conj:1}{\rm\cite[Conjecture 5.8]{Sun}} Let $m$ and $n$ be positive integers, and let $0\leqslant k<n$. Then
\begin{align}\label{eq:sun3}
\frac{[km+1]!}{([k]!)^m}\sum_{h=0}^{n-1}q^h {h\brack k}^m\equiv 0\pmod{[n]},
\end{align}
where $[n]!=[n][n-1]\cdots [1]$, and so
\begin{align}\label{eq:sun4}
\frac{(km+1)!}{(k!)^m}\sum_{h=0}^{n-1}{h\choose k}^m\equiv 0\pmod{n}.
\end{align}
\end{conj}

Conjecture \ref{conj:1} for $m=1$ is easy, and Sun himself is also able to prove the $m=3$ case of this conjecture.

The aim of this paper is to prove Conjecture \ref{conj:1} for arbitrary $m$ by establishing the following more general form.
\begin{thm}\label{thm:new}
Let $m,n\geqslant 1$, and let $a_1,\ldots,a_m$ be nonnegative integers. Then
\begin{align}\label{eq:new1}
\frac{[a_1+\cdots+a_m+1]!}{[a_1]!\ldots[a_m]!}\sum^{n-1}_{h=0}q^h\prod_{i=1}^m{h\brack a_i}
\equiv 0\pmod{[n]},
\end{align}
and so
\begin{align}\label{eq:new2}
\frac{(a_1+\cdots+a_m+1)!}{a_1!\cdots a_m!}\sum^{n-1}_{h=0}\prod_{i=1}^m{h\choose a_i} \equiv 0\pmod{n}.
\end{align}
\end{thm}

It is clear that, when $a_1=\cdots=a_m=k$, the congruences \eqref{eq:new1} and \eqref{eq:new2}
reduce to \eqref{eq:sun3} and \eqref{eq:sun4}, respectively.

For $m=2$, we shall prove the the following stronger result.
\begin{thm}\label{thm:new2}
Let $a$ and $b$ be nonnegative integers and $p$ a prime with
$p>\max\{a,b\}$. Then
\begin{align}
\frac{[a+b+1]!}{[a]![b]!}\sum_{h=0}^{p-1}q^h{h\brack a}{h\brack b}
\equiv (-1)^{a-b} q^{ab-{a\choose 2}-{b\choose 2}}[p]\pmod{[p]^2},  \label{eq:new2-1}
\end{align}
and so
\begin{align*}
\frac{(a+b+1)!}{a!b!}\sum_{h=0}^{p-1}{h\choose a}{h\choose b}\equiv
(-1)^{a-b}p \pmod{p^2}.
\end{align*}
\end{thm}

\section{Proof of Theorem \ref{thm:new} }
For any $m,n\geqslant 1$ and nonnegative integers $a_1,\ldots,a_m$, let
\begin{align*}
S_n(a_1,\ldots,a_m)
=\frac{[a_1+\cdots+a_m+1]!}{[n][a_1]!\ldots[a_m]!}\sum^{n-1}_{h=0}q^h\prod_{i=1}^m{h\brack a_i}.
\end{align*}
To prove \eqref{eq:new1}, it is equivalent to show that $S_n(a_1,\ldots,a_m)$ is a polynomial in $q$ with integer coefficients.
By \cite[(3.3.9)]{Andrews}, we have
\begin{align}
\sum_{h=0}^{n-1}q^h{h\brack a_1}={n\brack a_1+1}q^{a_1},
\end{align}
and so
\begin{align}
S_n(a_1)=\frac{[a_1+1]}{[n]}{n\brack
a_1+1}q^{a_1}={n-1\brack a_1}q^{a_1}.  \label{eq:sn-1}
\end{align}
This proves the $m=1$ case.

For $m\geqslant 2$, by the $q$-Chu-Vandermonde summation formula
(which is equivalent to \cite[(3.3.10)]{Andrews})
\begin{align}
\sum_{k=0}^n {a\brack k}{b\brack n-k}q^{k(b-n+k)}={a+b\brack n},
\end{align}
we have
\begin{align*}
{h\brack a_{m-1}}{h\brack a_{m}}
&={h\brack a_{m-1}}\sum_{k=0}^{a_m}{h-a_{m-1}\brack k}{a_{m-1}\brack a_{m}-k}q^{k(a_{m-1}-a_m+k)} \\
&=\sum_{k=0}^{a_m}{h\brack a_{m-1}+k}{a_{m-1}+k\brack a_m}{a_m\brack k}q^{k(a_{m-1}-a_m+k)}.
\end{align*}
It follows that
\begin{align}\label{eq:double-sum}
S_n(a_1,\ldots,a_m)
&=\frac{[a_1+\cdots+a_m+1]!}{[n][a_1]!\ldots[a_m]!} \nonumber\\
&\quad{}\times\sum^{n-1}_{h=0}q^h\prod_{i=1}^{m-2}{h\brack a_i}
\sum_{k=0}^{a_m}{h\brack a_{m-1}+k}{a_{m-1}+k\brack a_m}{a_{m}\brack k}q^{k(a_{m-1}-a_m+k)}.
\end{align}
Exchanging the summation order in \eqref{eq:double-sum}, and noticing that
\begin{align*}
&\hskip -2mm
\frac{[a_1+\cdots+a_m+1]![a_{m-1}+k]!}{[a_1+\cdots+a_{m-1}+k+1]![a_{m-1}]![a_m]!}{a_m\brack k}
={a_1+\cdots+a_m+1\brack a_{m}-k}{a_{m-1}+k\brack a_{m-1}},
\end{align*}
 we obtain the following recurrence relation:
\begin{align}
S_n(a_1,\ldots,a_m)
&=\sum_{k=0}^{a_m}{a_1+\cdots+a_m+1\brack a_{m}-k}{a_{m-1}+k\brack a_{m}}{a_{m-1}+k\brack a_{m-1}}q^{k(a_{m-1}-a_m+k)}\nonumber\\
&\quad{}\times S_{n}(a_1,\ldots,a_{m-2},a_{m-1}+k).  \label{eq:sn-2}
\end{align}
The proof then follows easily by induction on $m$.

\section{Proof of Theorem \ref{thm:new2} }
By \eqref{eq:sn-2} and  \eqref{eq:sn-1}, we obtain
\begin{align*}
\frac{[a+b+1]!}{[a]![b]!}\sum_{h=0}^{p-1}q^h{h\brack a}{h\brack b}
=[p]\sum_{k=0}^{b}{a+b+1\brack b-k}{a+k\brack a}{a+k\brack b}{p-1\brack a+k}q^{k(a-b+k)+a+k}.
\end{align*}
Noticing that, if $0\leqslant a+k\leqslant p-1$, then
\begin{align*}
{p-1\brack a+k}=\prod_{i=1}^{a+k}\frac{1-q^{p-i}}{1-q^i}
\equiv \prod_{i=1}^{a+k}\frac{1-q^{-i}}{1-q^i}
=(-1)^{a+k} q^{-{a+k+1\choose 2}}  \pmod{[p]}.
\end{align*}
Moreover, since $p>\max\{a,b\}$, we have ${a+k\brack a}\equiv 0\pmod{[p]}$ if $a+k\geqslant p$ and $k\leqslant b$.
This means that, for $0\leqslant k\leqslant b$, we always have
\begin{align*}
{a+k\brack a}{p-1\brack a+k}
\equiv {a+k\brack a}(-1)^{a+k} q^{-{a+k+1\choose 2}}  \pmod{[p]}.
\end{align*}
Therefore, to prove \eqref{eq:new2-1}, it suffices to show that
\begin{align}
\sum_{k=0}^{b}{a+b+1\brack b-k}{a+k\brack a}{a+k\brack b}(-1)^k q^{k(a-b+k)+a+k-{a+k+1\choose 2}} \label{eq:new3}
=(-1)^b q^{ab-{a\choose 2}-{b\choose 2}},
\end{align}
which is just a special case of the $q$-Pfaff-Saalsch\"utz's identity (see \cite[3.3.12]{Andrews}):
\begin{align}
\sum_{k=0}^{n}\frac{(x;q)_k(y;q)_k(q^{-n})_k{q^k}}{(q;q)_k(z;q)_k(xyq^{1-n}/z;q)_k}=\frac{(z/x;q)_n(z/y;q)_n}{(z;q)_n(z/xy;q)_n},
\label{eq:qPfaff}
\end{align}
where $(a)_n=(1-a)(1-aq)\cdots (1-aq^{n-1})$.\\

In fact, replacing $(x,y,z,n,k)$ by $(q^x,q^x,q^{-a-x},a+b+1,k+1)$ in \eqref{eq:qPfaff}, we get
\begin{align}
&\hskip -2mm
\sum_{k=-1}^{a+b}\frac{(q^x;q)_{k+1}(q^x;q)_{k+1}(-1)^{k+1}}{(q^{-a-x};q)_{k+1}(q^{3x-b};q)_{k+1}}
{a+b+1\brack k+1}q^{{k+1\choose 2}-(k+1)(a+b)} \nonumber\\
&=\frac{(q^{-a-2x};q)_{a+b+1}(q^{-a-2x};q)_{a+b+1}}{(q^{-a-x};q)_{a+b+1}(q^{-a-3x};q)_{a+b+1}}. \label{eq:qPfaff-2}
\end{align}
It is easy to see that, for $k\geqslant 0$,
\begin{align*}
&\lim_{x\to 0}\frac{(q^x;q)_{k+1}}{(q^{-a-x};q)_{k+1}}=(-1)^{a+1}q^{a+1\choose 2}{k\brack a},\\[5pt]
&\lim_{x\to 0}\frac{(q^x;q)_{k+1}}{(q^{-b+3x};q)_{k+1}}=\frac{(-1)^b}{3} q^{b+1\choose 2}{k\brack b},\\[5pt]
&\lim_{x\to 0}\frac{(q^{-a-2x};q)_{a+b+1}(q^{-a-2x};q)_{a+b+1}}{(q^{-a-x};q)_{a+b+1}(q^{-a-3x};q)_{a+b+1}}=\frac{4}{3}.
\end{align*}
Letting $x\to 0$ in \eqref{eq:qPfaff-2}, we are led to
\begin{align*}
\sum_{k=a}^{a+b}{k\brack a}{k\brack b}{a+b+1\brack k+1}(-1)^k q^{{k+1\choose 2}+{a+1\choose 2}+{b+1\choose 2}-(k+1)(a+b)}=(-1)^{a-b},
\end{align*}
which is clearly equivalent to \eqref{eq:new3}.

\section{An open problem}
By Faulhaber's formula (see \cite{GZ,GRZ,Knuth}), it is not hard to see that, for positive integers $m$ and $n$, there holds
\begin{align*}
(2m+2)!\sum_{h=0}^{n-1}h^{2m+1}\equiv 0\pmod{n^2}.
\end{align*}
We end this paper with the following conjecture.
\begin{conj}Let $m,n$ and $k$ be positive integers with $m\geqslant k$. Then
\begin{align*}
\frac{((2k+1)(2m+1)+1)!}{((2k+1)!)^{2m+1}}\sum_{h=0}^{n-1}{h\choose 2k+1}^{2m+1}\equiv 0\pmod{n^2}.
\end{align*}
\end{conj}

\vskip 5mm \noindent{\bf Acknowledgments.} This work was partially
supported by the Fundamental Research Funds for the Central
Universities and the National Natural Science Foundation of China
(grant 11371144).

\end{document}